\documentclass[12pt,leqno]{article}

\usepackage{amsmath,amsfonts,amssymb,amsthm,amscd}
\date{}

\newcommand\Z{\mathbb{Z}}

\newcommand\A{\mathbb{A}}

\newcommand\C{\mathbb{C}}

\newcommand\N{\mathbb{N}}

\newcommand\R{\mathbb{R}}

\newcommand\Q{\mathbb{Q}}

\newcommand\qp{\mathbb{Q}_p}

\newcommand{\Zp}{\Z_p}

\newcommand\fa{\mathfrak{a}}

\newcommand\fb{\mathfrak{b}}

\newcommand\pp{\mathfrak{p}}

\newcommand\fn{\mathfrak{n}}
\newcommand\fm{\mathfrak{m}}
\newcommand\gog{\mathfrak{g}}

\newcommand\fk{\mathfrak{k}}

\newcommand\ffh{\mathfrak{h}}

\newcommand\Ff{\mathcal{F}}
\newcommand\Hh{\mathcal{H}}
\newcommand\Aa{\mathcal{A}}

\newcommand\CPp{\mathcal{P}}

\newcommand\Ll{\mathcal{L}}

\newcommand\Res{\mathrm{Res}}
\newcommand\ind{\mathrm{ind}}

\newcommand{\ba}{\backslash}
\newcommand\rg{\rightarrow}
\newcommand\lgr{\longrightarrow}

\newcommand\Hom{\mathrm{Hom}}

\newcommand{\Lie}{\mathrm{Lie}}

\newcommand{\nequiv}{\equiv\kern-4.5mm/}

\def\adots{\mathinner{\mkern1mu\raise1pt\vbox{\kern7pt\hbox{.}}
\mkern2mu\raise4pt\hbox{.}
\mkern2mu\raise7pt\hbox{.}\mkern1mu}}

\newtheorem{monlem}{\textsc{Lemma}}[section]

\newtheorem{maprop}[monlem]{\textsc{Proposition}}
\newtheorem{moncoro}[monlem]{\textsc{Corollary}}
\newtheorem{montheo}[monlem]{\textsc{Theorem}}

\numberwithin{equation}{section}

\begin{document}

\title{On the Eisenstein functoriality in cohomology for maximal parabolic subgroups}

\author{Laurent Clozel}

\maketitle

\section*{Abstract}
In his paper, 'On torsion in the cohomology of locally symmetric varieties', Peter Scholze has introduced a new, purely topological method to construct the cohomology classes on arithmetic quotients of symmetric spaces of rational reductive groups originating from the cohomology of the similar quotients of Levi subgroups of maximal parabolic subgroups. We extend this construction beyond the cases he considers, and, in the complex case, to the cohomology of local systems. 

\section{Introduction}

Let $G$ be a reductive group defined over $\Q$, $A=A_G$ the neutral component of the group of real points in a maximal split central torus of $G$, $K_\infty \subset G(\R)$ a maximal compact subgroup. We are interested in the cohomology of the quotients $\Gamma\ba X$, where $X=G(\R)/A_G K_\infty$ is the symmetric space, and $\Gamma\subset G(\R)$ is a congruence subgroup. For our purposes, it will be better to consider the ad\`elic version, i.e., the quotients
\begin{equation}
S_K = G(\Q)\ba G(\A)/A_GK_\infty K
\end{equation}
where $\A$ denotes the ad\`eles of $\Q$, $\A_f$ the finite ad\`eles and $K\subset G(\A_f)$ is a compact open subgroup.

The quotient (1.1) is a finite union of spaces $\Gamma\ba X$. It is known from the work of Borel and Franke that $H^\bullet(S_K,\C)$ can be computed by automorphic forms : in fact
$$
H^\bullet (S_K,\C)= H^\bullet(\gog,K_\infty;\Aa(G_K))
$$
where $\Aa$ denotes the space of automorphic forms on
$$
G_K= G(\Q)\ba G(\A)/A_GK)\,,
$$
and the $(\gog ,K_\infty)$ cohomology is defined  by  the action of $G(\R)$ by right translations. See \cite{BW},\cite{F}.

There is a decomposition, due to Langlands,
\begin{equation}
\Aa(G_K) = \bigoplus_{P} \Aa_P(G_K)
\end{equation}
where $P$ runs over the association classes of $\Q$--parabolic subgroups, and $\Aa_P(G_K)$ is the space obtained from Eisenstein series ``induced from $P$''.

When the space of cusp forms of $G$ is understood,  this yields the computation of the part relative to~$G$ :
$$
H^\bullet(\gog,K_\infty;\Aa_G(G_K)) = H^\bullet(\gog,K_\infty;\Aa_{cusp}(G_K)) = \bigoplus_\pi H^\bullet (\gog,K_\infty;\pi)
$$
where $\pi$ ranges over the cuspidal summands of $L^2(G_K)$ with non--trivial cohomology. Harder, and then Schwermer, have proposed a program aiming at constructing cohomology classes in $\Aa_P(G_K)$ for $P\neq G$, obtained as differential form--valued Eisenstein series `` induced'' from cusp forms on $M(\Q)\ba M(\A)$, $P=MN$ being the Levi decomposition. Important results have been proved by Harder \cite{Ha1, Ha2, Ha3, Ha4}, who in particular solved the problem completely for $GL(2)$, and by Schwermer \cite{Sch1, Sch2} and others. For a more extensive survey of the known results,  see e.g. Grbac \cite{Grbac}. However, a general construction of the ``Eisenstein classes'' does not exist.

In his paper, ``On torsion in the cohomology  of locally symmetric varieties'' \cite{Scho},  Scholze introduced a new method to construct the ``Eisenstein'' classes associated to  \textit{maximal} parabolic subgroups\footnote{Scholze does not seem to think that his method was new. It is new and probably yields new results: see \S 4.}, a topological method relying on the Borel--Serre compactification rather than Eisenstein series. Scholze limited himself to two groups, $Sp(g,F)$ for a totally real field $F$ and $U(F)$ for a quasi--split group $U$ of even rank $2n$ associated to a $CM$ quadratic extension $E/F$. (In these cases, he considered the maximal parabolic subgroups with Levi subgroups, respectively, $GL(g,F)$ and $GL(n,E)$.) Our purpose here is to show that the method is perfectly general when ``Eisenstein functoriality'' for a \textit{maximal} parabolic subgroup is considered.

We carry out Scholze's construction  in two cases : in characteristic $\ell$ (or more generally with an Artin ring of coefficients $k$) for cohomology with trivial coefficients, in \S~2. The results are Theorem~2.4, Theorem~2.5\footnote{We refer the reader to the text, as it would be tedious to introduce the relevant notation here.}. Then, for complex cohomology, in \S~3 where we deal (as is natural in this case) with local systems coming from an arbitrary (complex) representation of $G$. The main results are Theorem~3.4, Proposition~3.5, Theorem~3.6. They imply that there exist cohomology classes in $H^\bullet(S_K)$ associated to all classes in the inner cohomology $H_!^\bullet(S_{K_M}^M)$ for all Levi subgroups $M$ in maximal parabolic subgroups $P=MN$, $S_{K_M}^M$ being the associated ad\`elic quotient. (Proposition~3.5 does not seem to follow directly from the analytic construction of Eisenstein cohomology when it applies, since we can directly relate classes with $\bar{\Q}$--coefficients.)

In \S 4, we compare the results obtained with the (known) consequences of the Harder--Schwermer method. In particular, if $G=GL(n)$, we examine whether the eigenclasses ``on $G$'' originating from the \textbf{cuspidal} cohomology of $GL(m)\times GL(m) (n=2m)$ can be obtained by the formation of Eisenstein classes.

Finally, \S 5 explores the apparent limit of the method : if $P\subset G$ is not maximal,  we explain why the natural generalisation of Scholze's method fails. (It would be interesting to find a topological argument in this case.) I am indebted to Leslie Saper for providing the proof of a crucial property of the Borel--Serre compactification (Theorem 5.2.) Its proof will appear later.

In conclusion, I emphasise that this is only the natural development of an idea entirely due to Scholze. I thank him, and also G\"unter Harder, Lizhen Ji, Colette Moeglin, Benjamin Schraen, Leslie Saper, Jack Thorne and David Vogan, for useful communications.

\section{Eisenstein cohomology form maximal para\-bolic subgroups : Scholze's argument}

\setcounter{equation}{0}

\subsection{} In this section we will  extend to general reductive groups the proof of Scholze \cite[V.2]{Scho}
 showing the existence of ``Eisenstein'' cohomology classes (or rather, eigencharacters) originating from the inner cohomology of maximal parabolic subgroups.

The argument, topological, relies on the Borel--Serre compactification. We first fix some notations.

Le $G$ be a connected reductive group of positive rank over $\Q$. We denote by $A_G$ the (topological) neutral component of the set of $\R$--points of the maximal $\Q$--split torus in the center of $G$. Fix a compact--open subgroup $K$ in $G(\A_f)$. We will assume that $K$ is
decomposed : $K=\Pi K_p$. Following Newton  and Thorne \cite{NT}, we also assume $K$ \textbf{neat} in the strong sense of Pink \cite{P} : let $\rho$ be a faithful representation of $G$ over $\Q$. We fix an algebraic closure $\bar{\Q}$ of $\Q$ and an embedding $\bar{\Q} \rightarrow \bar{\Q}_p$ for all $p$. If $g=(g_p) \in G(\A_f)$, consider for each $p$ the torsion subgroup $\Gamma_p$ of the group generated in $\bar{\Q}_p$ by the eigenvalues of $g_p$. (Thus $\Gamma_p\subset \bar{\Q})$. Then $g$ is \textbf{neat} if $\bigcap \Gamma_p=\{1\}$. A compact open subgroup $K$ is neat if all its elements are neat. If $G$ is a connected linear algebraic group, $H$ a  subgroup, $M$ a quotient, of $G$, and $K\subset G(\A_f)$ is neat, then $K\cap H(\A_f)$ and the image of $K$ in $M(\A_f)$ are neat (ibid.). Such a subgroup is neat in the more usual sense $(K\cap G(\Q)$ has no element of finite order but $1$),  and they form a basis for the compact--open subgroups.

 Let $K_\infty$ be a maximal compact subgroup of $G(\R)$. We consider
\begin{eqnarray*}
S_K   &:= & G(\Q)\ba G(\A)/A_G K_\infty K\\
      &= &G(\Q) \ba (X_G \times G(\A_f))/K
\end{eqnarray*}
where $X_G=G(\R)/{A_GK_\infty}$ is the symmetric space of $G(\R)$. \footnote{More precisely, of $G(\R)/A_G$. Note that $A_G$ depends on the $\Q$- structure.} (Thus $S_K$ is a finite union 
\begin{equation}
S_K=\coprod \Gamma_i \ba X_G
\end{equation}
where the $\Gamma_i$ are congruence subgroups of $G(\Q)$ ; the quotients are smooth.)

Let $P=MN \subset G$ be a $\Q$--parabolic subgroup, with a Levi decomposition. There are similar quotients associated to $P$ and $M$ :
\begin{eqnarray*}
S_{K_P}^P &=& P(\Q)\ba P(\A) / A_M\ K_{\infty,M} K_P\,,\\
S_{K_M}^M &=&M(\Q) \ba M(\A) /A_M\ K_{\infty,M}\ K_M
\end{eqnarray*}
where $K_{\infty,M}\subset M(\R)$ is maximal compact, and $K_P$,  $K_M$ are compact--open in the sets of finite adelic points. If $K_P=P(\A_f) \cap K$ and $K_M$ is the projection of $K_p$, they are, as we saw, neat.

There is a natural projection
$$
\pi_P : S_{K_P}^P \rightarrow S_{K_M}^M
$$
which realises the quotient $S_{K_P}^P$ as a fiber bundle with compact fiber\break $N(\Q)\ba N(\A)/K_N$, where $K_N=K_P \cap N(\A_f)$. Note that the fiber is of the form $\Gamma_N\ba N(\R)$, $\Gamma_N\subset N(\Q)$ a congruence subgroup. The quotients $S_{K_P}^P$, $S_{K_M}^M$ have ``finite'' expressions as~in~(2.1).

We now consider the Borel--Serre compactification $S_K^{BS}$ of $S_K$. We start with the Borel--Serrre ``bordification'' of $X_G$. This is a disjoint union
$$
X_G^{BS}= \coprod_P e(P)
$$
where $P$ runs over the $\Q$--parabolic subgroups of $G$ ; $e(G)=X_G$. Let $\CPp$ be the set of $\Q$--parabolic subgroups, and $\CPp'$ the set of \textbf{proper} parabolic subgroups. Then $\coprod\limits_{P\in \CPp'} e(P)$ is the boundary, $\partial X_G^{BS}$, of the manifold with corners $X_G^{BS}$. For $P\in \CPp'$, $e(P)=X_G/A_P$ where we write $A_P$ for $A_M$, and the action of $A_P$ (``geodesic flow'') on $X_G$ commutes with the left action of $P(\R)$ \cite{BS}. We define
\begin{equation}
S_K^{BS} = G(\Q)\ba (X_G^{BS} \times G(\A_f))/K\,.
\end{equation}
Since $G(\A_f)/K$ is discrete, and $G(\Q)$ acts without fixed points, it is a  manifold with corners. We can give two descriptions of $S_K^{BS}$. First write
$$G(\A_f) = \coprod_g G(\Q) g\, K\ (g\in G(\A_f))\,,
$$
the finite union that leads to the expression (2.1). Then (2.2) easily yields
\begin{equation}
S_K^{BS} = \coprod_g \Gamma_g\ba X_G^{BS}
\end{equation}
where $\Gamma_g = G(\Q)\cap  gKg^{-1}$. Thus $S_K^{BS}$ is a finite union of ``classical'' Borel--Serre compactifications.  (In particular $S_K^{BS}$, with the natural topology coming from (2.2) and the  topology of $X_G^{BS}$, is compact.) On the other hand, the set of $\Q$--parabolic subgroups, modulo $G(\Q)$, is finite. Let $\mathcal{Q}$ be a set of representatives.  Then from (2.2) we get
\begin{equation}
S_K^{BS} = \coprod_{P\in \mathcal{Q}} G(\Q)\ba (e(P) \times G(\A_f)) /K\,.
\end{equation}
Each piece is the part of $S_K^{BS}$ contributed by the ``parabolic subgroups of type $P$''. It can be further decomposed. Denoting it by $S_K^{BS}(P)$, we remark that $P(\A_f)\ba G(\A_f)/K$ is finite since $P(\A_f)\ba G(\A_f)$ is the set of points of a projective variety with values in $\A_f$, with a transitive action of $G(\A_f)$. Let $\{h\}$ be a set of representatives. Then
\begin{equation}
S_K^{BS}(P) = \coprod_h P(\Q)\ba (e(P)\times P(\A_f))/{K_P(h)}\,,
\end{equation}
with $K_P(h) = P(\A_f) \cap h K h^{-1}$. In turn we have (by properties of the geodesic action)
\begin{eqnarray*}
e(P) &=& X_G/A_P\\
&=& (P(\R) K_{G,\infty}/K_{G,\infty}) /A_P
\end{eqnarray*}
($K_{G,\infty}$ being suitably chosen with respect to $P$, as in the definition of the geodesic action)
$$
= P(\R)/A_M K_{M,\infty}\,,
$$
so each piece of (2.5) is of the form
$$
P(\Q)\ba P(\A) /A_M \ K_{M,\infty}K_P = S_{K_P}^P\,,
$$
with $K_P=K_P(h)$. Thus $S_{K_P}^P$ occurs as one of the pieces of the compactification. We note that (2.4) - in particular when passing to the limit for $K\rg 1$ - allows one to see the cohomology of $S_K^{BS}(P)$ - rather,  the limit - as an induced representation of $G(\A_f)$. We will use a simple aspect of this. 

Recall that $\coprod\limits_{P\in \mathcal{Q}'} S_K^{BS}(P)$, where $\mathcal{Q}'$ is the set of proper representatives, is the boundary of $S_K^{BS}$. Its dimension is $\dim(X_G)-1$. If $P$ is a $\Q$--parabolic subgroup with split component $A_P=A_M$, the dimension of $e(P)$ is $\dim(X_G) -(\dim A_P-\dim A_G)$.  In particular, the quotients $S_K^{BS}(P)$ for $P$ a \textbf{maximal} $\Q$--parabolic  subgroup are the open cells of $\partial S_K^{BS}$.

\subsection{} We now consider the cohomology of these spaces. Although the first arguments could be given for any ring of coefficients, we will assume that $k$ is a local Artin ring. (We are mostly interested in the case where $k$ is a field ;  this however allows us to consider, as does Scholze, the case of $k=\Z/\ell^n\Z$ for a prime $\ell$.)

Let $S$ be the (finite) set of primes such that $K_p \subset G(\Q_p)$ is not hyperspecial. We  consider the Hecke algebras
\begin{eqnarray*}
\Hh^S(G) &=& C_c (K^S\ba G(\A_f^S)/K^S,\Z)\\
&=& \bigotimes\limits_{p\notin S}C_c(K_p\ba G(\Q_p)/K_p,\Z)
\end{eqnarray*}
where $K^S = \prod\limits_{p\notin S}K_p$, and the tensor product is restricted. Similarly we have $\Hh^S(P)$ and $\Hh^S(M)$.

We can assume that our set $\mathcal{Q}$ of representatives of the parabolic subgroups is the set of standard parabolic subgroups, i.e., those containing a fixed, minimal parabolic subgroup $P_0$. Then there is a choice of maximal compact subgroups $K_p^0$ (all primes) such~that
$$
G(\Q_p) = K_p^0P(\Q_p) = P(\Q_p)K_p^0
$$
(all $p$, all $P\in \mathcal{Q}$). We assume $K_p=K_p^0$ so chosen for $p\notin S$. In this situation those are natural $\Zp$--structures on all groups $G$, $P$, $M$, $N\dots$ for $p\notin S$, and we assume that the measures on all groups $X$ give mass 1 to $X(\Zp)$. We have natural maps
$$
\rho : \Hh^S(G) \rightarrow \Hh^S(P)
$$
given by the restriction of functions, and
$$
\lambda : \Hh^S(P) \rightarrow \Hh^S(M)
$$
given by 
$$
\varphi \mapsto \int_{N(\Q_p)} \varphi(mn)dn:=\lambda\varphi(m)\,.
$$
They are morphisms of algebras. Note that the  ``constant term'' $\lambda$ is not normalised.

These algebras act naturally on cohomology spaces. For instance, a double coset $K_p g K_p$ in $\Hh_p(G)$ acts on $H^i(S_K^{BS},k)$ by the correspondence 
\begin{equation}
\begin{matrix}
S^{BS}(K\cap g K g^{-1}) &\xrightarrow[R_g]{} &S^{BS}(g^{-1}Kg\cap g K)\\
\downarrow &&\downarrow\\
S^{BS}(K)&&S^{BS}(K)
\end{matrix}
\end{equation}
$R_g$ being right translation by $g\in G(\qp)$, using the description (2.2). Moreover, these correspondences preserve the decomposition (2.4).

In particular we see that $\Hh^S(G)$ acts on $H^\bullet(\partial\, S_K^{BS},k)$, compatibly with the map $H^\bullet(S_K^{BS},k) \rg H^\bullet(\partial\, S_K^{BS},k)$.

Assume now that $P\in \CPp$ is maximal. Then $S_K^{BS}(P)$, and in particular its component $S_{K_{P}}^P$, is open in $\partial\,S_K^{BS}$. Thus we obtain a map $j_* : H_c^i(S_{K_P}^P,k)\rg H^i(\partial\, S_K^{BS},k)$. We also have the projection $\pi : S_{K_P}^P\rg S_{K_M}^M$, whence $\pi^* : H_c^i(S_{K_M}^M,k) \rg H_c^i(S_{K_P}^P,k)$ since the fibres of $\pi$ are compact.

\vskip2mm

{\monlem{

$(i)$ The map $j_*$ is equivariant under $\Hh^S(G)$, $\Hh^S(G)$ acting naturally on $H^i(\partial\, S_K^{BS},k)$ and by composition with $\rho : \Hh^S(G) \rg \Hh^S(P)$ on $H^i(S_{K_P}^P,k)$.

$(ii)$ The map $\pi^*$ is equivariant under $\Hh^S(P)$, $\Hh^S(P)$ acting on $H_c^i(S_{K_M}^M,k)$ by composition with~$\lambda$.}}

\vskip2mm

\textbf{Remark}.--- This is due to Scholze \cite[Lemma V.2.3]{Scho}, who does not, however, give a proof. I am greatly indebted to Newton and Thorne for providing this proof.

\vskip2mm

For the first assertion, consider the open injection $S_{K_P}^P \hookrightarrow \partial \, S_K^{BS}$. We must consider the action of a double coset $K_p\,g\, K_p$, by  the correspondence (2.6). However we have assumed that, for $p\notin S$, $G(\qp)=K_pP(\qp)= P(\qp)K_p$. We can therefore realise the correspondence (2.6) on $\partial S_K^{BS} \subset S_K^{BS}$ by taking $g\in P(\qp)$. But the action on (the cohomology of) $S_{K_P}^P$ is then given by the same correspondence, and this is compatible with the map between the cohomology spaces.

The second assertion is more difficult, but it is the statement proved by Newton and Thorne \cite[Cor. 3.9]{NT}. The comparison with their paper calls for the following comments; we refer to their paper for some notations, which we do not need here.

(1) Newton and Thorne prove a stronger statement, concerning the actions of $\Hh^S(P)$ and  $\Hh^S(M)$ on $R\Gamma_c(S^P_{K_P},k)$,  $R\Gamma(S^P_{K_P},k),$ $R\Gamma_c(S^M_{K_M},k)$, $R\Gamma_(S^M_{K_M},k)$ seen as elements of suitable derived categories. Our statement is simply obtained by taking instead the cohomology of the relevant complexes.

(2) In our case the modules $B$ (a $K_S$-module) and $A$ (a $K_{M,S}$-module) figuring in their statement coincide with the trivial $k$-module; the corresponding sheaves are the constant $k$-sheaves.

(3) Their definition of the action of the Hecke algebras on the cohomology spaces is given in their \S 2 (see their Proposition 2.10 and the subsequent comment.) It coincides with the definition we have given, by correspondences (see their Lemma 2.19.)

Consider now the diagram
\begin{equation*}
\begin{matrix}
H_c^i (S_{K_P}^P,k) &\lgr &H^i(\partial\, S_K^{BS},k) &\lgr &H^i(S_{K_P}^P, k)\,.\\
&j_*	&&	j^*
\end{matrix}
\end{equation*}
Set $j_c = j^*j_* : H_c^i (S_{K_P}^P,k) \rg H^i(S_{K_P}^P,k)$. Its image is by definition the compactly supported part of the cohomology, or 'inner' cohomology, denoted by  $H_!^i(S_{K_P},k)$. The diagram yields naturally a surjective map
\begin{equation}
Im (j_*) \rightarrow\!\!\!\!\rightarrow H_!^i(S_{K_P}^P,k)\,.
\end{equation}
We now give the rest of the arguments, again due to Scholze. We consider the full diagram
\begin{equation}
\begin{matrix}
H_c^i(S_{K_P}^P,k) &\xrightarrow[j_c\ ]{} &H^i(S_{K_P}^P,k)\\
\uparrow\pi^*	&&	\downarrow \pi_*\\
\noalign{\vskip2mm}
H_c^i(S_{K_M}^M,k) &\xrightarrow[j_c^M]{} &H^i(S_{K_M}^M,k)
\end{matrix}
\end{equation}

A priori $H^\bullet S_{K_P}^P,k) = H^\bullet(S_{K_M}^M,R\pi_*k)$. 
This is associated to the local system $H^\bullet(\Ff)$ where $\Ff \cong \Gamma_N\ba N(\R)$ is the fiber of $\pi$. However $\Ff$ is covered by the contractible space $\tilde{\Ff}=N(\R)$, and therefore we get a map $H^\bullet(\Ff) \rg H^\bullet(\tilde{\Ff})=k$ of local systems on $S_{K_M}^M$, whence a section $H^\bullet(\Ff) \rg H^0(\Ff)=k$. Thus $\pi_*$ sends $H^i(S_{K_P}^P,k)$ to $H^i(S_{K_M}^M,k)$. The diagram (2.8) clearly commutes.

\vskip2mm

{\monlem $\pi_*$ is equivariant, $\Hh(P)$ acting naturally on $H^i(S_{K_P}^P,k)$ and via $\lambda$ on $H^i(S_{K_M}^M,k)$.}

\vskip2mm

This is again \cite[Cor. 3.9]{NT}.

Finally, we obtain the following corollary. The commutativity of the diagram (2.7) implies that the map $\pi_* : H_!^i (S_{K_P}^P,k) \rg H_!^i(S_{K_M}^M,k)$ is surjective. We finally have :

\vskip2mm

{\moncoro There exists a surjective map
$$
Im\, j_* \lgr H_!^i(S_{K_M}^M,k)\,,
$$
where $j_* : H_c^i(S_{K_P}^P,k) \rg H^i( \partial\, S_K^{BS},k)$, equivariant for the action of $\Hh^S(G)$ on the  left--hand side, and its action via $\lambda\circ \rho$ on the right--hand side.
}
\vskip2mm

Finally, consider the long exact sequence in cohomology for the manifold with boundary $S_K^{BS}$ :
\begin{equation}
\begin{matrix}
\cdots \lgr H_c^i(S_K,k) \lgr H^i(S_K^{BS},k) \lgr H^i(\partial\, S_K^{BS},k)\lgr\\ \\
 \lgr H_c^{i+1}(S_K,k) \lgr \cdots
 \end{matrix}
\end{equation}

The submodule $j_*\ H_c^i(S_{K_P}^P,k) := H$ of $H^i(\partial\, S_K^{BS},k)$ admits a filtration
$$
0 \lgr H' \lgr H \lgr H'' \lgr 0
$$
with $H'\subset H^i(S_K^{BS},k) = H^i (S_K,k)$ and $H'' \subset H_c^{i+1}(S_K,k)$. We finally have  by Cor.~2.3 :

\vskip2mm

{\montheo Each irreducible subquotient of  the $\Hh_G^S$--module $H_!^i(S_{K,M}^M,k)$ is a subquotient of $H_c^{i+1}(S_K,k)$ or of $H^i(S_K,k)$.
}

\vskip2mm

We recall that the action of $\Hh_G^S$ on $H_!^i(S_{K_M}^M,k)$ is through the map $\lambda\circ \rho$. In the ``classical'' case $(k=\C)$, this is the map associated to non--normalized induction from $M(\A^S)$ to $G(\A^S)$ (through $P$).

Note that all the cohomology spaces are finite over $k$. Let $\fm$ be the maximal ideal of $k$, and $\kappa=k/\fm$. All modules $H$ considered can be filtered~:
$$
0 = \fm^r \, H \subset \fm^{r-1} \, H \subset \cdots \subset \fm H \subset H\,,
$$
the quotients being vector spaces over $\kappa$, and preserved by the Hecke algebras. On each quotient $\Hh_G^S$ acts via the quotient $\Hh_G^S \otimes \kappa$. Denote by $\fm_{\Hh}$ a maximal ideal of $\Hh_G^S \otimes \kappa$. We deduce :

\vskip2mm

{\montheo Assume $H_!^i(S_{K_M}^M,k)_{\fm_{\Hh}}\not= 0$. Then $H_c^{i+1}(S_K,k)_{\fm_{\Hh}}\not= 0$ or $H^i(S_K,k)_{\fm_{\Hh}}\not= 0$.}

\vskip2mm

 For a more precise result, relating the images of $\Hh_M$ and $\Hh_G$ in $End (H^i(S^M_{K_M},k)$ and $End (H^i(S_K,k)$, see \cite[Cor.5.2.4]{Scho}.

\section{The complex case : coefficient systems}

\subsection{} We now return to the constructions of \S 2.2, but we consider the case where $k=\C$, introducing coefficient systems. We assume given an algebraic, complex representation of $G$, $L$ ; it defines naturally a local system $\Ll$ on $S_K$. Since $S_K \subset S_K^{BS}$ is a homotopy equivalence, $\Ll$ extends to $S_K^{BS}$, and then restricts to $S_K^{BS}(P)$ for any $P$. Recall (2.5) that $S_K^{BS}(P)$ is a union of components of the form
$$
\begin{matrix}
P(\Q)\ba e(P) \times P(\A_f)/K_P(h)\,,\\ 
\noalign{\vskip2mm}
e(P)= \ P(\R)/A_M K_{M,\infty}
\end{matrix}
$$
which  in turn are a union of quotients
$$
\Gamma_P\ba P(\R)/A_M K_{M,\infty}
$$
for congruence subgroups  $\Gamma_P$. The quotient $P(\R)/A_MK_{M,\infty} = N(\R) X(M)$ is simply connected.

The component $e(P)$ of $X^{BS}$ can be retracted to $X$ (see e.g. Schwermer \cite[1.9]{Sch1} and then we see that $\Ll|_{S_K^{BS}(P)}$ is given on each component by the restriction of $L$ to $\Gamma_P$. We will simply denote by $\Ll$ the local system obtained on each of the spaces $S_K$, $S_K^{BS}$, $S_K^{BS}(P)$\dots

Returning to \S 2.2, we first recall that the diagram of maps (2.6) again yields naturally a map $H^i(S^{BS}(K),\Ll) \rg H^i(S^{BS}(K),\Ll)$. For simplicity we describe this in the case of $S_K$. The total space of the vector bundle $\Ll$ on $S_K$ is
\begin{equation}
G(\Q)\ba (X\times (G(\A_f)/K) \times L\,,
\end{equation}
$G(\Q)$ acting diagonally on the three factors. The only non--obvious map in (2.6) is the effect of $R_g$. Write $K'=K\cap gKg^{-1}$, $K''=g^{-1}Kg\cap K$.
We need a map $L_y \rg L_{R_gy}$, for a point $y\in S(K')$. A representative of $L_y$ in the quotient (3.1) is $(x,h) \times L$ where $(x,h)\in S\times G(\A_f)$ is a representative of $y$. The map $R_g$ is $(x,h) \mapsto (x,hg)$. Similarly $L_{R_gy}$ is represented by $(x,hg)\times L$. The map $(x,h,\ell) \mapsto (x,hg,\ell)$ $(\ell\in L)$ descends to the quotient (3.1). (This is of course well-known ; see Harder \cite{Ha5} for a fuller description of the action of the Hecke algebra.)

As before Lemma 2.1, we obtain compatible actions of $\Hh^S(G)$ on\break $H^\bullet(\partial \, S_K^{BS},\Ll)$ and $H^\bullet(S_K^{BS},\Ll)$.  For $P$ maximal, we obtain again a map $j_* : H_c^i(S_{K_P}^P,\Ll) \rg H^i(\partial\, S_K^{BS},\Ll)$.

The first part of Lemma 2.1 is proved as before :

\vskip2mm

{\monlem The map $j_*$ is equivariant under $\Hh^S(G)$, acting on $H_c^i(S_{K_P}^P,\Ll)$ by composition with $\rho$.} 

\vskip2mm

However, part (ii) of Lemma 2.1 is not sufficient in this case, as $\Ll$ (when non trivial) is not defined on $S_{K_M}^M$. In general, Harder and Schwermer have described $H^\bullet(S_{K_P}^P,\Ll)$ as a cohomology space for $S_{K_M}^M$. This is given by the degeneracy of the Leray spectral sequence for the fibration $S_{K_P}^P \rg S_{K_M}^M$ by compact nilmanifolds. Let $\fn$ be the real Lie algebra of $N$. For each $j$, $H^j(\fn,L)$ is an $M(\R)$--module and defines a local system $\Hh^j(\fn,L)$ on $S_K^M$. Then
\begin{equation}
H^i(S_{K_P}^P,\Ll)= \bigoplus_{j+k=i} H^k(S_{K_M}^M\,, \Hh^j(\fn,L))\,.
\end{equation}
This is due to Harder \cite{Ha1} and Schwermer \cite{Sch1}. Note that $\Hh^S(M)$ acts on the left--hand side. Harder and Schwermer prove the isomorphism (3.2) for classical quotients, i.e., for the components of our ad\`elic quotients.

However these classical quotients (respectively for $P$ and $M$) are in bijection : use strong approximation for the nilpotent group~$N$.

\vskip2mm

{\monlem $(i)$ The isomorphism $(3.2)$ is true for compactly supported cohomology on both sides.

$(ii)$ It is equivariant under $\Hh^S(P)$, acting on the right via~$\lambda$.}

\vskip2mm

Par $(i)$ is clear since it is the (degenerate) Leray spectral sequence, and the fibres are compact. For $(ii)$, we have to consider the proof  of Harder and Schwermer.

Let $\fm^0$ denote the Lie algebra of $M^0$, where $M(\R)=M^0A_M$ is the Langlands decomposition. Let $\pp^0=\fm^0\oplus \fn$. Let $\frak{k}_M$ be the Lie algebra of $K_{M,\infty}=K_{P,\infty}$. As described by Schwermer \cite[p. 49-50]{Sch1}, there is, for each $i$, $j$, a natural injection
\begin{equation}
\begin{matrix}
\eta^i : \Hom_{K_{M,\infty}}(\Lambda^j(\fm^0/\fk_M)\,,\ C^\infty(M(\Q)\ba M(\A)/A_MK_M) \otimes H^i(\fn,L))\\
\noalign{\vskip2mm}
\lgr \Hom_{K_{M,\infty}}(\Lambda^{i+j}(\pp^0/\fk_M)\,,\ C^\infty(P(\Q)\ba P(\A)/A_MK_P)\otimes L)\,.
\end{matrix}
\end{equation}

In fact Schwermer proves this statement for classical quotients $\Gamma_M\ba M(\R)/A_M$, $\Gamma_P\ba P(\R)/A_M$. Again, this implies the ad\`elic variant.

The map $\eta^i$ is obtained from the injection 
$$
\Lambda^\bullet  (\fm^0/\fk_M)^* \otimes \Lambda^\bullet (\fn)^* \lgr \Lambda^\bullet (\pp^0/k_M)^*
$$
of dual spaces, and from the natural injection
$$
C^\infty (M(\Q)\ba M(\A)/A_M K_M) \lgr C^\infty (P(\Q)\ba P(\A)/A_MK_P)
$$
given by the projection $P\rg M$. It induces an isomorphism in cohomology, yielding~(3.2).

We note that the same injection obtains between the spaces of compactly supported, smooth functions. The proof of (3.3) (\cite[Theorem 2.7]{Sch1}) now yields the same result for cohomology with compact support. The spaces of smooth functions on the ad\`ele groups in (3.2) receive, respectively, an action of $\Hh^S(M)$ and $\Hh^S(P)$ by right translations, which then yield the actions on the cohomology of the respective quotients. Therefore the proof of $(ii)$ is completed by the following easy lemma.

\vskip2mm

{\monlem Denote by $I$ the natural injection \;  $C_c^\infty(M(\Q)\ba M(\A)/A_MK_M) \rg C_c^\infty(P(\Q)\ba P(\A)/A_M K_P)$. Then, if $\varphi\in \Hh^S(P)$ and $f\in C_c^\infty(M(\Q)\ba M(\A)\break /A_M K_M)$ :
$$
I(\lambda(\varphi)f) = \varphi\ I(f))
$$
}

\vskip2mm

We can omit the compact--open subgroups and work with the functions on the full ad\`elic quotients. Denote by $pr$ the projection $P\rg M$ (for these quotients). Then, for $g\in P(\A)$ :
\begin{equation*}
\begin{matrix}
\varphi I(f)(g) &=&\displaystyle\int_{P(\A_f)} I(f)(gg')\varphi(g')dg'\hfill\\
\noalign{\vskip2mm}
&=&\displaystyle\int_{P(\A_f)} f(pr(gg'))\varphi (g')dg'\,,
\end{matrix}
\end{equation*}
while for $m\in M(\A)$ :
\begin{equation*}
\begin{matrix}
\lambda(\varphi) f(m) &=&\displaystyle\int_{M(\A_f)} f(mm')(\lambda\varphi)(m')dm'\hfill\\
\noalign{\vskip2mm}
&=&\displaystyle\int_{M(\A_f)N(\A_f)} f(mm')\varphi (m'n')dm'dn'\,,\hfill\\
I(\lambda(\varphi)f)(g) &=&\displaystyle\int_{M(\A_f)N(\A_f)}f(pr(g)m') \varphi(m'n')dm'dn'\,.
\end{matrix}
\end{equation*}
Since $pr(gg')=pr(g)pr(g')=pr(g)m'$, the property is clear.

\subsection{} As in \S 2, the diagram
$$
H_c^i (S_{K_P}^P,\Ll) \xrightarrow[j_*]{} H^i(\partial\, S_K^{BS},\Ll) \xrightarrow[j^*]{} H^i(S_{K_P}^P,\Ll)
$$
yields a surjective map
$$
Im(j_*) \lgr\!\!\!\!\rg H_!^i(S_{K_P}^P,\Ll)\,.
$$
of $\Hh^S(G)$--modules. However, since (3.2) is true both for $H^i$ and $H_c^i$, the right--hand side~is
\begin{equation}
\bigoplus_{j+k=i} H_!^k (S_{K_M}^M\,,\Hh^j(\fn,L))\,.
\end{equation}
We can now consider the long exact sequence for $S_K^{BS}$ as in (2.9), to conclude that the subquotients of (3.4), on which $\Hh^S(G)$ acts via $\lambda \circ \rho$, occur in $H^i(S_K,\Ll)$ or $H_c^{i+1}(S_K,\Ll)$.

Rather than by localising at a maximal ideal, we express the result, as usual in the complex case, in terms of characters of the (commutative) Hecke algebra $\Hh^S(G,\C)=\Hh^S(G)\otimes \C$.

Let $G_K$ denote the ad\`elic quotient of the group,
$$
G_K=G(\Q)\ba G(\A)/A_GK\,,
$$
so $S_K= G_K/K_{G,\infty}$. We then have an inclusion of representations of $G(\R)$
$$
L_{cusp}^2(G_K) \subset L_{dis}^2(G_K)\,,
$$
the space of cusp--forms inside the discrete spectrum, and consequently for the automorphic forms :
$$
\Aa_{cusp}(G_K) \subset \Aa_{dis}^2 (G_K)\,.
$$
The $(\gog,K_\infty)$--cohomology of these spaces gives the cuspidal and the $L^2$--cohomology (in the discrete sense : cohomology represented by $L^2$ harmonic forms) of $S_K$ ; moreover it is known~that
$$
H_{cusp}^\bullet(S_K) \subset H_!^\bullet(S_K) \subset H_{(2)}^\bullet(S_K)\,
$$
where $H^\bullet_{(2)}$ is the $L^2$-cohomology. Cf. Schwermer \cite{Sch3}. With coefficients $L$, we have 
$$
H_{cusp}^\bullet(S_K,\Ll) = H^\bullet(\gog,K\,;\ \Aa_{cusp}(G_K)\otimes L)\,.
$$
The same applies to $M$. Finally, the full cohomology of $S_K$ is, by Franke's result \cite{F} :
$$
H^\bullet (\gog, K_\infty\,;\ \Aa(G_K))
$$
for the full space of automorphic forms. The abstract form of our result is then :

\vskip2mm

{\montheo $(i)$ Assume $\chi$ is a character of $\Hh^S(M)$ occurring non--trivially in $H_!^k(S_{K_m}^M\,,\ \Hh^j(\fn,L))$ for some values of $k,j$ ; if $i=k+j$, $\chi \circ (\lambda\circ \rho):=\chi'$, a character of $\Hh^S(G)$, occurs non--trivially in $H^i(S_K,\Ll)$ or $H_c^{i+1}(S_K,\Ll)$.

$(ii)$ In particular, if $\chi$ occurs in $H_{cusp}^k(S_{K_M}^M,\, \Hh^j(\fn,L)=H^k(\fm,K_M;\break \Aa_{cusp}(M_{K_M}\otimes L')$ where $L'=H^j(\fn,L)$, $\chi'$ occurs in $H^i(S_K,\Ll)$ or $H_c^{i+1}(S_K,\Ll)$.}

\vskip2mm

(Note that by Poincar\'e dualit\'e $H_c^{i+1}
(S_K,\Ll)$ is also described by automorphic forms.)

The proof implies more : not only do the eigenspaces for $\Hh^S(M)$ (or $\Hh^S(P)$) yield eigenspaces for $\Hh^S(G)$, but the dimensions are conserved. Precisely, we have the following result. Denote by $M_\chi$ a generalised eigenspace associated to a character~$\chi$.

\vskip2mm

{\maprop Let $\chi$ be a character of $\Hh^S(M)$, and $\chi'$ the associated character of $\Hh^S(G)$. Then
$$
\dim\, H^i(S_K,\Ll)_{\chi'}
 + \dim\, H_c^{i+1} (S_K,\Ll)_{\chi'} \ge \sum_{i=j+k} \dim\, H_!^k (S_K^M,\, \Hh^j(\fn,L))_\chi
 $$ } 

\subsection{} Finally, we give the more explicit formulation of the result given by Kostant's theorem, following again Harder and Schwermer. Assume $L$ is irreducible, with highest weight $\nu\in \ffh_\C^*$, where $\ffh_\C$ is a Cartan subalgebra of $\gog_\C=\Lie (G(\R)/A)\otimes \C$. Let $\fm_\C = \Lie(M(\R)/A)\otimes \C$ and 
$\pp_\C = \Lie(P(\R)/A)\otimes \C$. We fix a Borel subalgebra $\fb_\C \subset \pp_\C$ ; $\pp_\C=\fm_\C\oplus \fn_\C$. We choose a split real form $\gog_r$ of $\gog_\C$, whence $\frak{h}_r$, $\fb_r$, $\pp_r=\fm_r+\fn_r$. (Note that $\fn_r=\fn=\Lie\, N(\R)$.) Let $W=W_G$, $W_M$ be the complex Weyl groups ; we assume given on $\frak{h}_r$ a scalar product invariant by $W$. Let $R_G$, $R_G^+$, $\Delta_G$ be the roots, positive roots, simple roots for $G$ ; similar notations for $M$. Let $\rho$ be the half--sum of roots for $R_G^+$. Finally, let $W^P$ be the set of Kostant representatives for $W_M\ba W$ :
$$
W^P = \{w\in W \mid w^{-1} \alpha \in R_G^+ \ \forall \alpha\in \Delta_M\},
$$
cf. \cite[\S 2.3]{Sch1}. We write $L_M^\xi$ for an irreducible complex representation of $M$ with highest weight $\xi$. Then
\begin{equation}
H^j(\fn,L) = \bigoplus_{w\in W^P\atop \ell(w)=j} L_M^{w(\nu+\rho)-\rho}
\end{equation}
(Kostant). We will use this description to show that (up to twists, see below) all the compactly supported cohomology of $S_{K_M}^M$ gives rise to cohomology of $S_K$ with corresponding action of the Hecke algebras.

Assume therefore that $\xi\in \ffh_r^*$ is a positive weight for $M$. (Note that we assume that all weights vanish on $\fa=\Lie \,A\subset \ffh_r$. The reader may as well suppose the split component $A$ trivial.) Following (3.5) we want to write
\begin{equation}
\xi
 +\rho = w(\nu+\rho)
 \end{equation}
with $\nu$ dominant for $G$ and $w\in W^P$. Since $\nu+\rho$ is strictly positive, this implies that $\xi+\rho$ must be regular for the roots of~$G$.

Assume this first. We can certainly write
$$
\xi+\rho=w\,\nu'
$$
with $w\in W$ and $\nu'$ positive (dominant) for $G$, and regular. Then for $\alpha\in \Delta_M$
$$
<\xi+\rho,\alpha> = <w\nu',\alpha> = <\nu',w^{-1}\alpha>
$$
is strictly positive, so $w^{-1} \alpha\in R_G^+$, and $w\in W^P$. Furthermore $\nu'$, being regular, verifies $<\nu',\alpha>\geq 1$ for $\alpha\in \Delta_G$, so $\nu'=\nu+\rho$.

In this case we can realise $L_M^\xi$ as a summand of (3.5), with this choice of $\nu$. Setting $j=\ell(w)$, and considering Theorem (3.5 $(ii)$, we assume that $H_!^k(S_{K_M}^M\,,\, \Ll_M^\xi)\not=0$. Then, with $i=j+k$, the corresponding character of $\Hh^S(G)\otimes \C$ will occur in $H^i$ or $H_c^{i+1}$ for~$S_K$.

Consider now the case where $\xi+\rho$ is not $G$--regular. For simplicity assume $A=A_G=\{1\}$. We certainly have $<\xi+\rho,\alpha>$ $>0$ for $\alpha\in \Delta_M$. The dual of the $1$--dimensional space $\fa_M \subset \ffh_r$, $\fa_M^*$, is identified with $X^*(M)\otimes \R$ where $X^*(M)\cong \Z$ is the group of $\Q$--rational characters of $M$. In particular, a character $\mu$ defines an element $T\in \fa_M^* \subset \ffh_r^*$. The representation $L_M^\xi\otimes \mu$ then has highest weight $\xi+\rho+T$. Recall that we have used a scalar product to identify $\ffh_r$ with its dual. Let $\beta\in \Delta_G$ be the unique root not in $\Delta_M$. Then, if we choose $T$ in a suitable half--line, $<\beta,T>>0$ and therefore $\xi+\rho+T$ is regular. Consequently, the twisted representations $L_M^\xi \otimes \mu$ (and the corresponding local system) verify our previous conditions.

If $\pi=\pi_\infty \otimes \pi_f$ is  a  representation  of  $M(\A)$ such that

$H^\bullet(\fm,K_M; \, \pi_\infty \otimes L_M^\xi)\not=0$, and   occurring in  the  discrete  spectrum, $\pi\otimes \vert {\mu} \vert^{-1}$ has cohomology with coefficients in $L_M^\xi\otimes \mu$ ; obviously the compactly supported cohomology classes correspond. Finally, we can strengthen Theorem~3.4 :

\vskip2mm

{\montheo Assume $H_!^k(S_{K_M}^M,\Ll_M)\not=0$ for any coefficient system, associated to $L_M$. Let $\chi:\Hh^S(M)\otimes \C \rg \C$ occur in this space. Then the character $\chi':\Hh^S(G)\otimes \C$ associated to an Abelian, unramified twist of $\chi$ occurs in $H^i(S_K,\Ll)$ or $H_c^{i+1}(S_K,\Ll)$, $\Ll$ being constructed as above.}

\subsection{} We now strengthen the results of \S 3.2, by controlling the rationality of the cohomology classes. Note that the local system $\Ll$, and all the cohomology spaces, are defined over $\overline{\Q}$. Precisely, we consider $\Ll$ as a local system of $\overline{\Q}$-vector spaces. Lemma~3.2 is clearly still true over $\overline{\Q}$. (The local system $\Hh^j(\fn,L)$ is a local system of $\overline{\Q}$--vector spaces.) We can now prove the following ``rational'' result. Denote, as before, by $j:S_{K_P}^P\rg \partial S_K^{BS}$ the natural injection, and let
$$
\delta : H^i (\partial S_K^{BS},\Ll) \lgr H_c^{i+1}(S_K,\Ll)
$$ 
be the connecting homomorphism in (2.9). Finally, denote by $\Res : H^i(S_K^{BS},\Ll)\break \rg H^i(\partial S_K^{BS},\Ll)$ the natural restriction.

Let $\chi:\Hh^S(G,\overline{\Q})\rg \overline{\Q}$ be a character, and, for each (finite--dimensional) $\Hh^S(G,\overline{\Q})$--module $M$, let $M_\chi$ be the generalised eigenspace associated to $\chi$. Finally, recall that $H^i(S_{K_P}^P,\Ll)=\bigoplus\limits_{j+k=i}H^k(S_{K_M}^M,\Hh^j(\fn,L))$.

\vskip2mm

{\montheo Let $\alpha\in H_!^i(S_{K_P}^P,\Ll)$ be a non--zero generalised eigenclass associated to $\chi$, for the map $\Hh^S(G)\rg \Hh^S(M)$. Then one of the following properties is true :

$(i)$ $\alpha=j^*\ \Res\,\beta$

\noindent for a (non--zero) class $\beta\in H^i(S_K^{BS},\Ll)_\chi = H^i(S_K,\Ll)_\chi$

$(ii)$ $\alpha=j^* \gamma$

\noindent for a class $\gamma\in H^i(\partial S_K,\Ll)_\chi$ such that
$$
\beta= \delta(\gamma) \in H_c^{i+1}(S_K,\Ll)_\chi\not=0\,.
$$}

\vskip2mm

We repeat that the cohomology is now taken with $\overline{\Q}$-coefficients. 

\section{Relation with Eisenstein cohomology}

\subsection{} We now check, at least is a simple case, the compatibility of this construction with the theory of Eisenstein cohomology.

We will consider only the simplest case (and the only perfectly general one) where "Eisenstein classes'' have been constructed : this is Schwermer's Theorem~4.11 in \cite{Sch1}. (Note that Schwermer's result is not limited to  maximal parabolic subgroups.)

Let $L=L^\nu$, the coefficient system for $G$, be given, and consider a cuspidal cohomology class $[\varphi]\in H_{cusp}^\bullet(S_{K_P}^P,\Ll)$. Assume it is associated to a representation $\pi$ of $M^0(\R)$ occurring in the cusp forms, where $M(\R)=M^0(\R)A_M$ is the Langlands decomposition, and, following (3.4) and (3.5), to $w\in W^P$ --- a ``class of type $(\pi,w)$'' in the terminology of Schwermer \cite[p. 82]{Sch1}. Schwermer constructs a (differential--form valued) Eisenstein series $E(\varphi,\lambda)$, where $\varphi$ is a harmonic representative of $[\varphi]$, and $\lambda\in \fa_M^* \otimes \C$.  \footnote{In fact $(\fa_M/\fa_G)^*\otimes \C$ ; similarly, consider $\fa_M/\fa_G$ in (4.1). All our linear forms are trivial on $\fa_G$.} The form $E(\varphi,\lambda)$ is constructed, by summation over $\Gamma_P\ba \Gamma$, from a form $\varphi_{\lambda} $ on $e'(P) \times A_P$: see Schwermer \cite[\S 3.3]{Sch1} If it is holomorphic at the point
\begin{equation}
\lambda
_0=-w(\nu+\rho)|_{\fa_M}\,,
\end{equation}
this yields a closed form on $S_K$, which represents a non--trivial cohomology class for $G$. Note that since $\pi$ is cuspidal unitary, the domain of convergence for $E(\varphi,\lambda)$ is given \cite[p. 86]{MW}~by\begin{equation}
<\mathrm{Re}  \,\lambda,\alpha^\vee>\ > \\ <\rho_G-\rho_M,\alpha^\vee>
\end{equation}
for each coroot $\alpha^\vee$, $\alpha\in \Delta_G$, occurring in $N_P$. (cf. \cite[ \S 6.3]{Sch1}).

Now consider the cohomology classes constructed in \S 3. We start, as in Theorem 3.5, with a character $\chi$ of $\Hh^S(M)\otimes \C$ occurring in the cohomology $H_!^\bullet(S_{K_M}^M,\Ll_M)$. Then (perhaps after an unramified, Abelian twist) it occurs in $H^\bullet(S_K,\Ll)$ for some coefficient system $L$. If $L_M$ is one of the representations $H^j(\fn,L)$ for some $j$ and $L$, the twist in not required.

\subsection{} We now limit ourselves to the case of $GL(n)$\footnote{For a detailed study of certain Eisenstein classes when all the Eisenstein series are holomorphic, in the case of $GL(n)$, see Harder-Raghuram \cite{HR}}. Consider the decomposition~(1.1)
\begin{equation}
\Aa(G_K)= \bigoplus_P \Aa_P(G_K)
\end{equation}
and the corresponding decomposition of the cohomology. The contribution of $\Aa_P(G_K)$ to $H^\bullet(S_K)$ is given by the Eisenstein series (and their residues) coming from cuspidal representations of $M(\A)$. Consider specifically the cuspidal cohomology $H_{cusp}^i(S_{K_M}^M,\C)$. It injects in $H_!^i(S_{K_M}^M,\C)$. In particular, if a character $\chi$ of $\Hh^S(M)$ occurs in $H_{cusp}^i(S_K^M,\C)$, the associated character $\chi'$ of $\Hh^S(G)$ occurs in $H^i(S_K)$ or $H_c^{i+1}(S_K)$ --- in the second case, the character $\tilde{\chi}'$ associated to the dual representation occurs in $H^{d-i-1}(S_K)$.

Assume now $G=GL(n)$, so $P$ has Levi subgroup $GL(n_1)\times GL(n_2)$. By the results of Jacquet and Shalika on strong multiplicity one, $\chi'$ (or $\tilde{\chi}'$) must occur in the cohomology given by the summand of (4.3) associated to~$P$.

We can try to construct it using Schwermer's theorem. Thus let $\omega=\omega_1\otimes \omega_2$ be a harmonic $i$--form on $S_{K_M}^M$, a product of the arithmetic quotients for $GL(n_1)$ and $ GL(n_2)$, that represents our cohomology class.  As recalled above (with different notation), we deduce from $\omega$ and $s \in \frak{a}^*_M \otimes \C$ a form $\omega_s$ on $e'(P) \times A_P$ and (in the appropriate range) an Eisenstein class $E(\omega, s)$. Since the cohomology with trivial coefficients for $M$ corresponds to $w=1$ in (4.1), we must evaluate the Eisenstein series $E(\omega_1\otimes \omega_2,s)$ at the point $s=-\rho|_{\fa_M}$.

Now assume further that the level is everywhere unramified. We also assume that $n_1=n_2= \frac{n}{2}$. The holomorphy of $E(\omega_1\otimes \omega_2,s)$ is governed by its constant term. This is given by
\begin{equation}
E_P(\omega_1 \otimes \omega_2,s) = (\omega_1 \otimes \omega_2)_s +(M(s)(\omega_1 \otimes \omega_2))_{\tilde{s}}
\end{equation}

where $\tilde{s}=(s_2,s_1)$ if $s=(s_1,s_2)$. The operator $M(s)$ can be decomposed as $M(s)=N(s) L(s)$ where $L(s$) is, according to Langlands, equal to the scalar

$$
\frac{L(\pi_1 \otimes \tilde{\pi}_2,s_1-s_2)}{L(\pi_1\otimes \tilde{\pi}_2,s_1-s_2+1)}. 
$$
Here the $L$--function is the Rankin $L$--function, and $\omega_1$, $\omega_2$ belong to the cuspidal representations $\pi_1$, $\pi_2$.  The terms $E_P(\omega_1 \otimes \omega_2,s), (\omega_1 \otimes \omega_2)_s $, can be considered as belonging to $\ind_{P(\A)}^{G(\A)} ((\pi_1\otimes \pi_2) \vert \det\vert^s)$, and the last term in (4.4) to  $\ind_{P(\A)}^{G(\A)} ((\pi_1\otimes \pi_2) \vert \det \vert^{\tilde{s}})$. See Arthur \cite{Ar} and Schwermer \cite [\S 4.2]{Sch1}

 Now $2\rho|_{\fa_M}=2\rho_N$ (where $P=MN$) is equal to $|\det\, m_1|^{n_1}|\det\, m_2|^{-n_2}$, so $-\rho|_{\fa_M}$ corresponds to $s_1-s_2=-\frac{n}{2}$. Thus the non--trivial part of the constant term is normalised by the $L$-function
$$
\frac{L(\pi_1\otimes \tilde{\pi}_2,-\frac{n}{2})}{L(\pi_1\otimes\tilde{\pi}_2,-\frac{n}{2}+1)} =
\frac{\varepsilon(\pi_1\otimes \tilde{\pi}_2)-\frac{n}{2})L(\pi_1\otimes\tilde{\pi}_2,1+\frac{n}{2})}{\varepsilon(\pi_1\otimes \tilde{\pi}_2,1-\frac{n}{2})L(\pi_1\otimes\tilde{\pi}_2,\frac{n}{2})}\,.
$$
This is always holomorphic (including for $n=2$), and non-zero if $n>2$, which we now assume. Moreover, we see that the evaluation  at $-\rho$ (for the usual parametrization of the Eisenstein series, corresponding to unitary induction) corresponds to the unnormalized constant term : $\Hh_S(G) \rightarrow \Hh_S(M)$, so this is compatible with our construction. 

 Since we are evaluating $M(s)$ on unramified functions, the normalised operator $N(s)$ is reduced to its Archimedean factor
 $$
 N_\infty(s): I(s):=\ind_P^G(\pi_1[s_1] \otimes (\pi_2[s_2]) \mapsto \ \ind_P^G(\pi_1[s_2] \otimes (\pi_2[s_2])
 $$
 ($\pi_1, \pi_2$ now denote the Archimedean components; $\pi[t] = \pi \otimes  \vert det  \vert ^t$) where induction is normalised.) We will simply denote this operator by $N(s)$.We are evaluating at $s_0= (-\frac{m}{2},\frac{m}{2})$; $N(s)$ is meromorphic, and Schwermer's construction relies on the holomorphy of $N(s)$ at $s_0$ - or more precisely, on the vector in 
$ \ind_P^G(\pi_1[-\frac{m}{2}] \otimes (\pi_2[\frac{m}{2}]) = I(s_0)$ deduced from $\omega_1Ê\otimes \omega_2$. 

Assume, for definiteness, that $m$ is even. The representation $I(s)$ is very explicit.  We have $\pi_1= \pi_2 := \pi$ and this representation is the tempered representation of $GL(m)$ having  non-trivial cohomology with trivial coefficients. Cf.\cite [\S 3.5]{Cl}. For $t \in \frac{1}{2} \N, t>0$, let $ \delta_t$ be the discrete series representation of $GL(2,\R)$ with Langlands parameter given on $\C^{\times}= W_{\C}$ by
$$
z \mapsto ((z/\bar{z})^t, (z/\bar{z})^{-t}). 
$$

Then 
$$
\pi= \delta_{ \frac{1}{2} } \times \delta_{ \frac{3}{2}}  \times \dots \times \delta_{\frac{m-1}{2}} ,
$$

 a tempered representation. (Here we use $\times$ to denote induction by blocks.) Thus $I(s)$ is
 $$
 ( \delta_{ \frac{1}{2} } \times \delta_{ \frac{3}{2}}  \times \dots \times \delta_{\frac{m-1}{2}}) [s_1] \times ( \delta_{ \frac{1}{2} } \times \delta_{ \frac{3}{2}}  \times \dots \times \delta_{\frac{m-1}{2}} )[s_2].
 $$
 
 For $\tilde{s_0} =(\frac{m}{2},- \frac{m}{2})$, we see that $I(\tilde{s_0})$ has a unique irreducible quotient, its Langlands quotient $Q:=Q(\tilde{s_0})$; $\pi$ being self-dual, $I(s_0)$ has a unique irreducible submodule. 
 
 Now $N(s) N(\tilde{s})=1;$ $N(s)$ is holomorphic at $\tilde{s}_0$; and  $N(s)$ is holomorphic as $s_0$ on $Q\subset I(s_0)$. Since the normalisation factor is holomorphic and non-zero, this is equivalent to the same property for the non-normalised operator, say $N'(s)$. The holomorphy of $N'(\tilde{s_0})$ is part of the construction of the Langlands quotient, which also implies that $Q$ is the unique irreducible submodule of $I(s_0)$; the holomorphy of $N(s)$  at $s_0$ on $Q\subset I(s_0)$ follows from the functional equation. Indeed, let $v\neq 0$ be an element of $Q \subset I(s_0).$ We can view $v$ as a vector independent of $s$ in the induced representation. We have $N(s) N(\tilde{s})v=v$. The vector  $w = N(\tilde{s})v$ is holomorphic and nonzero at $s_0$, for which value it is in $Q$. This implies that $N(s)$ is holomorphic at $s_0$ on $w$ and therefore on $Q$.

 However, $N(s)$ is \textit{not} holomorphic at $s_0$ on the full representation $I(s_0)$. This follows from \cite [I.2 Lemme]{MW2} : the holomorphy of $N(s)$ is equivalent to the irreducibility of $I(s_0)$; however, this representation is not irreducible as a consequence of a theorem of Speh and Vogan \cite [Thm. 6.15]{SV}.\footnote{I leave it to the reader to unravel the definitions of Speh and Vogan for our data.} In fact, one has a more complete result:
 
 {\monlem   Assume $v \in  I(\tilde{s_0})$ does not belong to $Q$.  Then $N(s)v$ has a pole at $s_0$.}

 \vskip2mm
 
Since $N(\tilde{s})$ is holomorphic at  $\tilde{s_0}$ and the two induced representations are dual, it defines a holomorphic family of invariant hermitian forms on the constant space of the induced representation. By a result of Vogan \cite[Theorem 3.8]{Vo} this form vanishes, at $s_0$, on $ker N(\tilde{s_0})$. Suppose $N(s) v$ is holomorphic at $s_0$. Then $N(\tilde{s}) N(s) v =v$ is holomorphic. Thus $v=N(s_0) N(\tilde{s_0}) v$ belongs to $Q=Im N(s_0)$.

\vskip2mm

Thus the only vectors in the induced representation on which the Eisenstein series is holomorphic are the vectors in $Q$. The infinitesimal character of $I(s_0)$ is equal to that of the trivial representation, and thus all its subquotients are candidates at having non-trivial cohomology.  We do not know in which subquotients the form deduced from $\omega_1, \omega_2$ may occur; however, there seems to be no reason that it will always belong to $Q$. However \footnote{Vogan, personal communication} it seems difficult to compute the cohomology of $Q$; in particular we do not know if one so obtains a space of the dimension given by Proposition 3.5.  Thus, even in this simple case, it is possible that we have obtained classes which are not (directly) obtained by Schwermer's construction. One should also recall that Scholze's construction, in general, will succeed for all classes in the inner cohomology, which will not always be represented by cusp forms.

\section{Non--maximal parabolic subgroups}

\subsection{} In this section we explain why Scholze's construction seems -- without a new idea -- limited to the case of maximal parabolic subgroups.

So assume $P\subset G$ is an arbitrary parabolic subgroup, contained in a maximal one, $Q$. We consider again the union of faces $S_{K_P}^P \subset \partial S_K^{BS}$. Recall that this is a union of faces $e'(P)$, of the form $\Gamma_P \ba e(P)$. The argument will concern one face at a time, so we work classically rather than in the ad\`elic formulation. We  have embeddings
$$
e'(P) \hookrightarrow \overline{e'(Q)} \subset \partial S_\Gamma^{BS}
$$
where $S_\Gamma^{BS}$ is now a component of $S_K^{BS}$. See \cite[Proposition 9.4]{BS}. By ``smoothing'' (see the Appendix to \cite{BS}) we can see $\partial S_\Gamma^{BS}$ as a smooth variety ; $e'(P)$ is then a locally closed submanifold. In particular, $e'(P)$ and $\partial S_\Gamma^{BS}$ are orientable and satisfy Poincar\'e duality. Let $j$ denote the embedding $e'(P)\rg \partial S_\Gamma^{BS}$. There still exists a direct image morphism
$$
j_* : H_c^i(e'(P)) \lgr H^{i+c}(\partial\, S_\Gamma^{BS})
$$
where $c$ denotes the codimension. If $d$ is the dimension of the symmetric space, $\partial\, S_\Gamma^{BS}:=\partial$ has dimension $d-1$ and $e'(P)$ $d-a$, where $a=\dim(A_P/A_G)$. Thus $c=a-1$. The map $j_*$ is defined by Poincar\'e duality :
$$
(j_*\, \alpha,\beta)_\partial = (\alpha,j^*\beta)_{e'(P)}\,.
$$
We forego coefficient systems and consider cohomology with complex coefficients. Here $\alpha\in H_c^i(e'(P))$, $\beta\in H^{d-a-i}(\partial)=H^{\dim(\partial)-c}(\partial)$. In order to imitate the previous argument, we would need to consider
$$
H_c^i(e'(P)) \xrightarrow[j_*]{} H^{i+c}(\partial)\xrightarrow[j^*]{} H^{i+c}(e'(P))\,.
$$
However :

\vskip2mm

{\maprop  Assume $P$ is not maximal. Then $j^*j_*=0$.}

\vskip2mm

This will follow from the following result, kindly communicated by Leslie Saper. Assume $P=P_1 \subset P_2 \subset P_3 \dots \subset P_a =Q$ is a sequence of parabolic subgroups, with $a_{i +1}= a_i-1, a_i = dim( A_{P_i}/A_G)$, and $Q$ maximal. Consider the associated embeddings $e'(P_i) \subset \overline{e'(P_{i+1})}.$

\vskip2mm

{\montheo \textnormal{(L. Saper)}.- $(i)$ The embedding $e'(P_i) \rightarrow \overline{e'(P_{i+1})} \subset \partial$ is $C^{\infty}$-homotopic to a smooth immersion $e'(P_i) \rightarrow e'(P_{i+1}). $

  $(ii)$ Consider the composed immersion $j:e'(P)\rg \partial$. Then $j$ is homotopic within $\partial$, and in fact within $e'(Q)$, to a map (in fact a smooth immersion) $k$ such that $Im(j) \cap Im(k)=\emptyset$. If $\omega \subset Im\, j$ is compact, we can even assume that the closure of $Im(k)$ does not meet $\omega$.}

\vskip2mm

We will construct $j_*\alpha$ that is ``supported on $Im\,j$'',  implying that $j^*(j_*\alpha)=k^*(j_*\alpha)=0$. We can proceed as follows.

Let $\tilde{\alpha}$ be a closed, compactly supported form on $e'(P)$ representing $\alpha$. Assume its support is contained in an open subset $U\subset e'(P)$ with compact closure. Since $e'(P) \mapsto \partial$ is obtained from a sequence of immersions in codimension $1$, the normal bundle to $e'(P)$ in $\partial$ is trivial. Therefore there exists a neighbourhood $V$ of $U$ in $\partial$, with compact closure, and a diffeomorphism $U\times I^c \xrightarrow[]{\sim} V$, $I=]-1,1[$. We can assume that $V$ does not meet $Im\,(k)$. Let $(y_r)$, $r=1,\ldots c$, be the coordinates on $I^c$. Consider the current on $V$ :
$$
\tilde{\gamma} = \delta_0(y)\tilde{\alpha} \wedge dy_1 \wedge \cdots \wedge dy_c
$$
where $\int_{I^c} \delta(y) \varphi(y)dy\equiv \varphi(0)$. It is then  easy to see that for any closed form $\tilde{\beta}$ on $\partial$, of degree $d-a-i$,
$$
\int_\partial \tilde{\gamma} \wedge \tilde{\beta} = \int_{e'(P)} \tilde{\alpha} \wedge j^*\tilde{\beta}\,.
$$

We can approximate $\tilde{\gamma}$, as a closed current with compact support in $V$, by closed forms $\tilde{\theta}$. We obtain cohomology classes arbitrarily close to $j_*\alpha$ ; $j_* H_c^i(e'(P))$ being finite--dimensional, we see that we can so obtain this whole space. Clearly the forms $\tilde{\theta}$ verify $k^*\tilde{\theta}=0$; since  $j$ and $k$ are homotopic this implies $j^*j_*=0$.

\textbf{Remark}.-- With a more thorough argument using currents, it may be possible to dispense with the second part of Theorem 5.2. However we think that Proposition 5.1 should remain true even with an arbitrary system of coefficients - at least a field $k$, as in $\S 2$, using Verdier duality. In this case the argument using the deformed embedding may be necessary. This is left to the reader.

\subsection{} Although this may be obvious, we remark that we cannot use our construction inductively to obtain cohomology for $G$ from the cohomology of a Levi subgroup. Indeed, we had to start with classes in $H_!^\bullet(S_{K_P}^P)$. The ``Eisenstein'' classes we constructed in $H^\bullet(S_K^G)$ come from the boundary, cf. (2.9) ; if they occur in $H^i(S_K^{BS})$ they do not belong to $H_!^i$ ; if they lie in $H_c^{i+1}(S_K)$ they are sent to $0$ in $H^{i+1}(S_K^{BS})= H^{i+1}(S_K)$. We cannot use ``induction by stages'' for a chain of parabolic subgroups !

\eject


\begin{thebibliography}{3}

\bibitem{Ar} Arthur, James Eisenstein series and the trace formula. Automorphic forms, representations and L-functions, pp. 253-274, Proc. Sympos. Pure Math., XXXIII,  Part 1, Amer. Math. Soc., Providence, R.I., 1979.
\bibitem{BS} Borel, A.; Serre, J.-P.  Corners and arithmetic groups. Comment. Math. Helv. 48 (1973), 436-491.
\bibitem{BW} Borel, Armand; Wallach, Nolan R. Continuous cohomology, discrete subgroups, and representations of reductive groups. Annals of Mathematics Studies, 94. Princeton University Press, Princeton, N.J., 1980.
\bibitem{Cl} Clozel, Laurent Motifs et formes automorphes: applications du principe de fonctorialit\'{e}.  Automorphic forms, Shimura varieties, and L-functions, Vol. I (Ann Arbor, MI, 1988), 77-159, Perspect. Math., 10, Academic Press, Boston, MA, 1990.
\bibitem{F} Franke, Jens Harmonic analysis in weighted L2-spaces. Ann. Sci. Ecole Norm. Sup. (4) 31 (1998), no. 2, 181-279.
\bibitem{Grbac}  Grbac, Neven Eisenstein cohomology and automorphic L-functions. Cohomology of arithmetic groups, 35-50, Springer Proc. Math. Stat., 245, Springer, Cham, 2018. 
\bibitem{Ha1} Harder, G. On the cohomology of discrete arithmetically defined groups. Discrete subgroups of Lie groups and applications to moduli (Internat. Colloq., Bombay, 1973), pp. 129-160. Oxford Univ. Press, Bombay, 1975.
\bibitem{Ha2}Harder, G. Eisenstein cohomology of arithmetic groups. The case GL(2). Invent. Math. 89 (1987), no. 1, 37-118. 
\bibitem{Ha3}Harder, G\"{u}nter  Some results on the Eisenstein cohomology of arithmetic subgroups of GLn. Cohomology of arithmetic groups and automorphic forms (Luminy-Marseille, 1989), 85-153, Lecture Notes in Math. 1447, Springer, Berlin, 1990
\bibitem{Ha4} Harder, G. Arithmetic aspects of rank one Eisenstein cohomology. Cycles, motives and Shimura varieties, 131-190, Tata Inst. Fund. Res. Stud. Math., 21, Tata Inst. Fund. Res., Mumbai, 2010.
\bibitem{Ha5} Harder, G. Cohomology of arithmetic groups, to appear.
\bibitem{HR} Harder, G. Raghuram, A. Eisenstein cohomology for $GL_N$ and the special values of Rankin-Selberg $L$-fuctions, Princeton University Press, annals of Math. Studies, 2020
\bibitem{La} Langlands, R. P. On the classification of irreducible representations of real algebraic groups. Representation theory and harmonic analysis on semisimple Lie groups, 101-170, Math. Surveys Monogr., 31, Amer. Math. Soc., Providence, RI, 1989.
\bibitem{MW} Moeglin, Colette; Waldspurger, Jean-Loup D\'{e}composition spectrale et s\'{e}ries d'Eisenstein. Une paraphrase de l'Ecriture. Progress in Mathematics, 113. Birkh\"{a}user Verlag, Basel, 1994. 
 \bibitem{MW2} Moeglin, C.; Waldspurger, J.-L. Le spectre r\'{e}siduel de GL(n). Ann. Sci. Ecole Norm. Sup. (4) 22 (1989), no. 4, 605-674.
\bibitem{NT} Newton, J., Thorne, J., Torsion Galois representations over CM fields and Hecke algebras in the derived category, preprint.
\bibitem{P}  Pink, Richard Arithmetical compactification of mixed Shimura varieties. Bonner Mathematische Schriften, 209. Universit\"{a}t Bonn, Mathematisches Institut, Bonn, 1990.
\bibitem{Scho}  Scholze, Peter On torsion in the cohomology of locally symmetric varieties. Ann. of Math. (2) 182 (2015), no. 3, 945-1066.  
\bibitem{Sch1}  Schwermer, Joachim Kohomologie arithmetisch definierter Gruppen und Eisensteinreihen. Lecture Notes in Mathematics, 988. Springer-Verlag, Berlin, 1983.

\bibitem{Sch2} Schwermer, Joachim Eisenstein series and cohomology of arithmetic groups: the generic case. Invent. Math. 116 (1994), no. 1-3, 481-511.
\bibitem{Sch3}  Schwermer, Joachim Cohomology of arithmetic groups, automorphic forms and L-functions. Cohomology of arithmetic groups and automorphic forms (Luminy-Marseille, 1989), 1-29, Lecture Notes in Math., 1447, Springer, Berlin, 1990.


\bibitem{SV} Speh, Birgit; Vogan, David A., Jr. Reducibility of generalized principal series representations. Acta Math. 145 (1980), no. 3-4, 227-299.
\bibitem{Vo} Vogan, David A., Jr. Unitarizability of certain series of representations. Ann. of Math. (2) 120 (1984), no. 1, 141-187.


\end{thebibliography}
\end{document}